\documentstyle{amsppt} 
\magnification 1200
\UseAMSsymbols
\hsize 5.5 true in
\vsize 8.5 true in
\parskip=\medskipamount
\NoBlackBoxes

\def\mathbb{\Bbb }

\def\mathcal{\Cal}

\def\supp{\text{\rm supp\,}}

\def\mod{\text{\rm {mod\,}}}
\def\ve{\varepsilon}
\def\vp{\varphi}

\def\snint{\raise2pt\hbox{$_{^\not}$}\kern-3.5 pt\int}

\def\snint{\raise2pt\hbox{$_{^\not}$}\kern-3.5 pt\int}

\TagsOnRight
\NoRunningHeads

\document
\parskip=\medskipamount
\topmatter
\title
On the Fourier-Walsh spectrum of the Moebius function
\endtitle
\author
J.~Bourgain
\endauthor
\address
Institute for Advanced Study, Princeton, NJ 08540
\endaddress
\email
bourgain\@ias.edu
\endemail
\abstract
We study the Fourier-Walsh spectrum $\{\widehat\mu(S); S\subset\{1, \ldots, n\}\}$ of the Moebius function $\mu$ restricted to $\{0, 1, 2, \ldots, 2^n -1\}\simeq 
\{0, 1\}^n$ and prove that it is not captured by levels $\{\widehat \mu(S)| \, |S|< n^{\frac 23-\ve}\}$.  An application to correlation with monotone Bolean functions is given.
\endabstract
\endtopmatter

\noindent
{\bf 0. Introduction}

This paper may be seen as a companion of \cite {B1} on the behavior of the Fourier-Walsh coefficients of the Moebius function $\mu$ 
restricted to a large interval $\{1, 2, \ldots, N\}, N=2^n$.
While in \cite {B1} we did establish nontrivial uniform upper bounds on the F-W coefficients of $\mu$, we are interested here in their distribution.
Our main result shows that $\mu$ cannot be captured by `low order' Walsh functions, more precisely

\proclaim
{Theorem 1} Let $\lambda>0$ be a fixed constant and $n_0\sim (\log n)^{-1-c\lambda} n^{2/3}$.
Then
$$
\sum_{|A|\leq n_0} |\widehat\mu (A)|^2 <(\log n)^{-\lambda}\tag 0.1
$$
where
$$
\widehat\mu (A) =\frac 1N\sum_{x\in \{0, 1\}^n} w_A(x) \mu\Big(\sum^{n-1}_{j=0} 2^jx_j\Big)\tag 0.2
$$
and
$$
w_A(x) =\prod_{j\in A} (1-2x_j).\tag 0.3
$$
\endproclaim

We note here that \cite{B1} does not provide a statement of this strength, at least not unconditionally.
Also to be mentioned is B.~Green's paper \cite{Gr} that contains a similar result, but with a slightly lower cutoff level $n_0$ (a small technical
upgrading of the argument in \cite{Gr} would allow to reach $n_0=n^{\frac 12-\ve}$).

The paper \cite{Gr}, \cite{B1} and also the present one are motivated by the general problem of understanding the computational complexity of $\mu$.
In particular, our exponent $\frac 23$ in Theorem 1 allows us to address a question posed by G.~Kalai on the (non)-correlation of $\mu$ with monotone
Boolean functions (\cite{B1} only provides a conditional proof of this fact).
We are invoking here an important result of Bshouty and Tamon \cite{B-T}, implying roughly that any monotone Boolean function on $\{0, 1\}^n$ 
has most of its Fourier-Walsh spectrum below level $0(\sqrt n)$.
Hence Theorem 1 implies

\proclaim
{Corollary 2}
The Moebius function does not correlate with monotone Boolean functions.
\endproclaim

Of course Theorem 1 is equally interesting in the context of circuit complexity and in particular B.~Green's result \cite {Gr} on $AC^o$-circuits.
The proof of Theorem 1 rests on the circle method and has similarities with the author's paper \cite{B2} on establishing a prime number with
prescribed binary digits $x_j=a_j\in \{0, 1\}$ for $j\in A\subset\{1, \ldots, n-1\}$.
In both cases, the main interest is braking the $\frac 12$-barrier (in \cite{B2}, we allow sets $A$ satisfying $|A|< n^{\frac 47-\ve}$, while a
previous result in \cite {H-K} is limited to $|A|<n^{\frac 12-\ve}$ for $A$ in general position).
Braking this $\frac 12$-barrier, unconditionally, seems to require a certain refinement of the `classical' method.
We should indeed refer to the results of Balog and Perelli \cite{B-P} that immediately imply a stronger version of Theorem 1 with $n_0$ replaced by
$n_1\sim \frac n{(\log n)^2}$, for some constant $c>0$, provided one assumes that there are no Siegel zero's (see \cite{B-P}, Theorem 3).
While we use part of the analysis in \cite{B-P}, additional work is required to estimate the effect of a possible Siegel zero in our problem (mostly
carried out in \S4 of the paper).
Significant here is the observation made in \cite{Gr} that the modulus $r$ of a primitive Siegel character is not a power of 2.
This feature is responsible for an amazing difference between the usual Fourier and Fourier-Walsh expansions of the Moebius function
regarding the effect of possible Siegel zero's.

We also want to point out that, while \cite{B2} is conceptionally close to the present paper, there are significant differences between
the problems studied and we were forced to reproduce several modified arguments from \cite{B2}.

Finally, the exponent $\frac 23$ in Theorem 1 may be a restriction of technical nature.
Our primary goal was to go beyond $\frac 12$ and more work may lead to further improvement.

\noindent
{\bf 1. Preliminaries}

Let $N=2^n, n\in\Bbb Z_+$ large and restrict the Moebius function $\mu$ to $[0, N[$.

Identifying $[0, N[$ with the Boolean cube $\Omega = \{0, 1\}^n$ by binary expansion
$$
x=\sum_{0\leq j<n} x_j 2^j\tag 1.1
$$
the Walsh system $\big\{ w_A; A\subset \{0, 1, \ldots, n-1\}\big\}$ is defined by $w_\phi =1$ and
$$
w_A(x) =\prod_{j\in A} (1-2x_j)=e^{i\pi\sum_{j\in A} x_j}.\tag 1.2
$$

The Walsh functions on $\Omega$ form an orthonormal basis \big(the character group of $(\mathbb Z/ 2\mathbb Z)^n$\big) and given a function $f$
on $\Omega$, we write
$$
f =\sum_{A\subset \{0, \ldots, n-1\}} \widehat f(A) w_A\tag 1.3
$$
where $\widehat f(A) =2^{-n}\sum_{n\in\Omega} f(n) w_A(n)$ are the Fourier-Walsh coefficients of $f$.

Understanding the size and distribution of those coefficients is well-known to be important to various issues, in particular in complexity theory
and computer science.
Roughly speaking, a F-W spectrum which is `spread out' indicates a high level of complexity for the function $f$.

Specifying $f=\mu|_\Omega$, we proved in \cite{B1} the uniform estimate
$$
|\widehat\mu (A)|< 2^{-n^{\frac 1{10}}} \ \text { for all } \ A\subset\{0, 1, \ldots, n-1\}.\tag 1.4
$$
We consider here the related but slightly different problem of how large $n_0\in\Bbb Z_+$ can satisfy
$$
\sum_{|A|\leq n_0} |\widehat\mu (A)|^2 =o(1).\tag 1.5
$$
Thus (1.5) means that $\mu$ is not captured by Walsh functions of weight below $n_0$.

For $0\leq\rho\leq 1$, define $K_\rho$ on $\{0,1 \}^n$ by
$$
K_\rho(x) =\prod^{n-1}_{0} (1+\rho -2\rho x_j).\tag 1.6
$$
Clearly $K_\rho\geq 0$ and $\int_\Omega K_\rho=1$.
Denote $T_\rho$ the convolution operator
$$
T_\rho f (x) =\int_\Omega f(x+y- 2xy) K_\rho(y)dy=\sum \widehat f(A) \rho^{|A|} w_A(x)\tag 1.7
$$
which is a contraction.

Assume $n_0\leq n$ and
$$
\sum_{|A|\leq n_0} |\widehat f(A)|^2 >c> n^{-C}.\tag 1.8
$$
Take $\rho = 1-\frac 1{n_0}$.
Hence, from (1.7)
$$
\langle f, T_\rho f\rangle > c\Big(1-\frac 1{n_0}\Big)^{n_0}> c'\tag 1.9
$$
and moreover
$$
\Vert T_\rho f\Vert_\infty \leq\Vert f\Vert_\infty\tag 1.10
$$
$$
\sum_{|A|>\ell} |\widehat{T_\rho f} (A)|^2 \leq \Big(1-\frac 1{n_0}\Big)^{2\ell}\Vert f\Vert^2_2.\tag 1.11
$$
Thus the advantage of $T_\rho f$ over $\sum_{|A|\leq n_0} \widehat f(A)w_A$ is to preserve (1.10).

We also need to involve the usual Fourier spectrum.

Let $h:\Bbb R\to \{1, -1\}$ be the 1-periodic function defined by
$$
\left\{
\aligned h=1 \, &\text { if } 0\leq x<\frac 12\\ h= -1 \,  &\text { if } \frac 12\leq x<1.\endaligned
\right.
$$
Hence
$$
w_A(x) =\prod_{j\in A} h\Big(\frac x{2^{j+1}}\Big).\tag 1.12
$$
Write for $x\in\Bbb Z$
$$
h\Big(\frac x{2^{j+1}}\Big) =\sum_r a_{r, j} \, e\Big(\frac {rx}{2^{j+1}}\Big) \ \text { where }
\sum_r |a_{r, j}|< Cj.\tag 1.13
$$
It follows that
$$
w_A(x)=\sum_r a_r e\Big(\frac {rx}{2^n}\Big) \ \text { with } \  \sum |a_r|< (Cn)^{|A|}.\tag 1.14
$$
Given $f$, the generating function
$$
S_f(\alpha )=\sum_{x=0}^{N-1} f(x) e(\alpha x)\tag 1.15
$$
is a 1-periodic function on $\Bbb R$.

By (1.14)
$$
|S_{w_A} (\alpha)|\leq \sum_r |a_r| \Big(\frac 1N+\Big\Vert\alpha+\frac r{2^n}\Big\Vert\Big)^{-1}
$$
and hence
$$
\Vert S_{w_A}\Vert_{L^1(\Bbb T)}< C\log N(Cn)^{|A|}.\tag 1.16
$$

In view of (1.11), (1.16), we deduce

\proclaim
{Lemma 1} Assume $f$ on $[0, N[$ satisfies $\Vert f\Vert_\infty\leq 1$.
Given $1<K<\frac n{2n_0}$ there is a decomposition
$$
T_\rho f= f_1+f_2\qquad \Big(\rho = 1-\frac 1{n_0}\Big) \tag 1.17
$$
such that
$$
\widehat {f_1} (A) =0 \ \text { for } \ |A|>K n_0\tag 1.18$'$
$$
and
$$
\Vert S_{f_1}\Vert_1 <(Cn)^{2Kn_0}\tag 1.18$''$
$$
and
$$
\Vert f_2\Vert_2 < e^{-K}  \sqrt N.
\tag 1.19
$$
\endproclaim

Let $f=\mu|_{[1, N[}$ and assume (1.8) holds for some $n_0=o(n)$.

Let $\rho=1-\frac 1{n_0}$ and denote $T_\rho[\mu]$ by $f$,
$$
S(\alpha)= \sum^N_1\mu(k) e(h\alpha) \ \text { and } \ S_f(\alpha) =\sum^N_1 f(k) e(k\alpha).\tag 1.20
$$

\bigskip

\noindent
{\bf 2. Minor Arcs Contribution}

We fix a parameter $B=B(n)$ which will be specified later.
At this point, let us just say that $\log B=o(n)$.

The major arcs are defined by
$$
\Cal M(q, a)= \Big[\Big|\alpha -\frac aq\Big|<\frac {B}{qN}\Big] \ \text { where \, $q<B$}.\tag 2.1
$$
Given $\alpha$, there is $q<\frac NB$ such that
$$
\Big|a-\frac aq\Big|< \frac B{qN} <\frac 1{q^2} \ \text { and \, $(a, q)=1$}.
$$
From Vinogradov's estimate (Theorem 13.9 in [I-K])
$$
\align
|S(\alpha)|&<\Big(q^{\frac 12} N^{-\frac 12}+q^{-\frac 12}+ N^{-\frac 15}\Big)^{\frac 12} (\log N)^4 N\\
&\ll \Big(\frac N{B^{\frac 14}}+\frac N{q^{\frac 14}}+N^{\frac 9{10}}\Big) (\log N)^4.\tag 2.2
\endalign
$$
Hence if $q\geq B$,
$$|S(\alpha)|\ll \frac N{B^{\frac 14}}(\log N)^4.\tag 2.3
$$
Denoting
$$
\Cal M=\bigcup_{q\leq B} \bigcup_{(a, q)=1} \Cal M(q, a)
$$
it follows that
$$
\Vert S\Vert_{L^\infty(\Bbb T\backslash \Cal M)} < CN(\log N)^4 B^{-1/4}.\tag 2.4
$$
From (1.9), (1.20)
$$
\int_0^1 S(\alpha)\overline{S_f}(\alpha) d\alpha>c'.\tag 2.5
$$
Fix $K$ and let $f=f_1+f_2$ be the decomposition from Lemma 1. By (1.18), (1.19), (2.4)
$$
\align
&\Big|\int^1_0S(\alpha) \overline{S_f} (\alpha)d\alpha- \int_{\Cal M }S(\alpha)\overline{S_f}(\alpha)d\alpha\Big|\leq\\
&\Vert S\Vert_2 .\Vert S_{f_2}\Vert_2 +\int_{\Bbb T\backslash \Cal M}|S| \ |S_{f_1}|\\
&< e^{-K}N+CN(\log N)^4 B^{-\frac 14}(Cn)^{2Kn_0}.
\endalign
$$
Taking $K\sim |\log c'|$ and
$$
\log B >C(\log n)n_0+C\log\log n\tag 2.6
$$
we obtain
$$
\Big|\int_{\Cal M} S(\alpha) \overline{S_f}(\alpha) d\alpha\Big|> \frac 12 c'.\tag 2.7
$$

\bigskip

\noindent
{\bf 3. Major Arcs Analysis (I)}

Consider the contribution
$$
\int_{\Cal M} S.\overline{S_f}.\tag 3.1
$$
At this stage we first recall some results from \cite{B-P}.
Following \cite{B-P}, let
$$
R=\exp \Big(c_1 \frac {\log N}{\log\log N}\Big)\tag 3.2
$$
with $c_1>0$ an appropriate constant.

There is an absolute constant $c_0>0$ such that for at most one primitive character $\Cal X(\mod q), q\leq R$, the Dirichlet $L$-function $L(s, \Cal
X)$ may vanish in the region
$$
\sigma> 1-\frac {c_0}{\log R} ; |t|\leq R^2.\tag 3.3
$$
If the exceptional character exists, then $L(s, \Cal X)$ has a unique zero $\beta$ in the region (3.3) and $\beta$ is simple and real.

Next, \cite{B-P} distinguishes the two cases.

\noindent
{\bf Case I.} The exceptional character exists and $\beta$ satisfies
$$
1-\beta<\frac {c_0}{2\log R}.\tag 3.4
$$

\noindent
{\bf Case II.} The exceptional character does not exist or, if it exists, $\beta$ satisfies
$$
1-\beta\geq \frac {c_0}{2\log R}.\tag 3.5
$$

\proclaim
{Lemma 2} (\cite{B-P}, Theorem 3).

If Case II, then
$$
\int_{\Cal M_1} |S(\alpha)|^2 d\alpha\ll _A\frac N{(\log N)^A}\tag 3.6
$$
with
$$
\Cal M_1 =\bigcup_{q\leq R} \bigcup_{(a, q)=1} \Big\{\alpha; \Big|\alpha-\frac aq\Big|< \frac {R(\log N)^{A+1}}{qN}\Big\}.\tag 3.7
$$
\endproclaim

Recalling the above definition of $\Cal M$, it follows from (3.6) that for
$$
B<R\tag 3.8
$$
$$
\int_{\Cal M} |S|.|S_f|\leq \Vert S_f\Vert_2.\Vert S|_{\Cal M_1}\Vert_2 <  C_AN(\log N)^{-A/2}\tag 3.9
$$
(contradiction). 
Hence, recalling (2.6), we obtain that
$$
\sum_{|A|<n_0} |\widehat\mu(A)|^2 =o(1) \ \text { for } \ n_0\sim \frac {\log N}{(\log\log N)^2}\tag 3.10
$$
conditional to the absence of Siegel zero's. 

The remainder of the paper aims at establishing an unconditional result.

We first perform the following manipulation of $f$.

Fix some
$$
\log B<m<\frac n{100}\tag 3.11
$$
(to be specified later) and partition $[1, n]$ in intervals $J_\alpha$ of size
$$
|J_\alpha|\sim m.\tag 3.12
$$
Let
$$
K_0=\Big(1+\frac {n_0m}n\Big) e^{3K}\tag 3.13
$$
and define for $A\subset\{0, 1, \ldots, n-1\}$
$$
\omega_\alpha(A)=
\left\{\aligned 1 \ &\text { if } \ |A\cap J_\alpha|\geq K_0\\ 0 &\text { \  otherwise.}\endaligned
\right.
$$
Let $f= f_1+f_2$ be the decomposition from Lemma 1 and estimate using (1.19) 
$$
\sum_A\omega_\alpha (A) |\widehat f(A)|^2 \leq e^{-2K}N+\sum_A \omega_\alpha (A) |\widehat f_1(A)|^2.\tag 3.14
$$
If $\widehat f_1(A) \not= 0$, then $|A|\leq K n_0$ and hence $\sum_\alpha \omega_\alpha(A) \lesssim \frac {Kn_0}{K_0}$, implying
$$
\sum_{\alpha\lesssim \frac nm}\Big[\sum_A \omega_\alpha (A) |\widehat f_1(A)|^2\Big] \lesssim \frac {Kn_0}{K_0} N.
$$
Therefore there is some $\alpha$ such that
$$
\frac n{4m} <\alpha < \frac n{2m}\ \text {  and } \ \sum_A\omega_\alpha(A) |\widehat f_1(A)|^2 \lesssim \frac {Kn_0m}{K_0n} N< e^{-2K} N.\tag 3.15
$$
The function
$$
g(x) =\sum_{|A\cap J_\alpha|< K_0} \widehat f(A) w_A\tag 3.16
$$
satisfies by (3.14), (3.16)
$$
\Vert f-g\Vert_2 < e^{-K} N\tag 3.17
$$
and also
$$
\Vert g\Vert_\infty < n^{K_0}.\tag 3.18
$$
From (2.7), (3.17),
$$
\Big|\int_{\Cal M} S(\alpha) \overline S_g (\alpha) d\alpha\Big|>\frac 13 c'.\tag 3.19
$$
We analyze
$$
\sum_{(a, q)=1} \int_{\Cal M(q, a)} S(\alpha) \overline S_g(\alpha) d\alpha\tag 3.20
$$
using multiplicative characters.

Expand
$$
e\Big(\frac {ak}q\Big)=\frac 1{\phi(q)} \sum_{\Cal X(\mod q)} \overline {\Cal X}(a) c_{\Cal X} (k)\tag 3.21
$$
where
$$
c_{\Cal X} (k) =\left\{\aligned &\overline{\Cal X}_1 \Big(\frac k{(k, q)}\Big) \frac {\vp(q)}{\vp(\frac q{(k, q)})} \mu\Big(\frac q{q_1 (k, q)}\Big)
\Cal X_1\Big(\frac q{q_1(k, q)}\Big) \tau (\Cal X_1) \ \text { if } \ q_1 \Big|\frac q{(q, k)}\\ &0 \ \text { otherwise}\endaligned\right.\tag 3.22
$$
where $\Cal X(\mod q)$ is induced by the primitive character $\Cal X_1(\mod q_1)$.

Let $0\leq \eta\leq 1$ be a suitable compactly supported smooth bumpfunction on $\Bbb R$ centered at $0$.
Writing $\alpha=\frac aq +\beta\in\Cal M(q, a)$ and using orthogonality, we obtain for (3.20)
$$
\frac B{qN} \frac 1{\phi(q)} \sum_{\Cal X(\mod q)} \sum_{k_1, k_2} c_{\Cal X} (k_1)\overline{c_{\Cal X}(k_2)} \, \mu(k_1)
\overline {g(k_2)}\, \widehat\eta \Big((k_1-k_2) \frac B{qN}\Big).\tag 3.23
$$
Assume the exceptional character $\Cal X^* (\mod r)$ exists, $r<B$, and denote $B$ the set of characters $\Cal X(\mod q), r|q$ that are induced by
$\Cal X^*$.

Decompose
$$
\spreadlines {8pt}
\align
&\int_{\Cal M} S. \overline S_g=\\
&\frac BN \sum_{q<B} \frac 1{q\phi(q)} \, \sum_{\Sb \Cal X(\mod q)\\ \Cal X\not\in B\endSb}\,  \sum_{k_1, k_2} c_{\Cal X} (k_1) 
\, \overline{c_{\Cal X}(k_2)} \, \mu(k_1)\, \overline {g(k_2)} \, \widehat\eta \Big((k_1-k_2)\frac B{qN}\Big)\tag 3.24\\
&+\frac BN \sum_{\Sb q< B\\ r|q\endSb} \, \frac 1{q\phi(q)} \sum_{k_1, k_2} c_{\Cal X}(k_1) c_{\Cal X} (k_2) \mu(k_1)\, \overline{g(k_2)}
\, \widehat\eta \Big((k_1-k_2)\frac B{qN}\Big)\tag 3.25\\
\endalign
$$
where in (3.25), $\Cal X$ is the unique character $(\mod q)$ induced by $\Cal X^*$.

Again by orthogonality and Cauchy-Schwarz inequality
$$
\align
&(3.24) =\frac BN \sum_{q<B} \frac 1q \, \sum_{(a, q)=1} \, \sum_{k_1, k_2} \frac 1{\phi (q)} \Big[\sum_{\Sb \Cal X\not\in \Cal B\\ (\mod q)\endSb}
\overline{\Cal X} (a) c_{\Cal X} (k_1)\Big] e\Big( -\frac {ak_2}q\Big) \mu(k_1)\overline{g(k_2)} \, \widehat \eta \Big((k_1-k_2)\frac B{qN}\Big)\\
&\leq \sum_{q<B}\, \sum_{(a, q)=1} \int_{|\beta|< \frac B{qN}}\Big|\sum^N_{k=1} \mu(k) \frac 1{\phi(q)} \Big[\sum_{\Sb \Cal X\not\in \Cal B\\ (\mod
q)\endSb} \overline{\Cal X}(a) c_{\Cal X}(k)\Big] e(k\beta)\Big| \ \Big|S_g \Big(\frac aq+\beta\Big)\Big| d\beta\\
& \leq \Big\{ \sum_{q<B} \, \sum_{(a, q)=1} \int_{|\beta|<\frac B{qN}}\Big|\sum^N_{k=1} \mu(k) \frac 1{\phi(q)} \Big[\sum_{\Sb \Cal X\not\in\Cal B\\
(\mod q)\endSb} \overline {\Cal X} (a) c_{\Cal X}(k)\Big] e(k\beta)\Big|^2 d\beta\Big\}^{\frac 12}.\Vert S_g\Vert_2\\
&<(3.26)^{\frac 12} \sqrt N.\tag 3.26
\endalign
$$

Restrict $q\sim Q<B$.
Performing the $\beta$-integration, we bound (3.26) by
$$
(\log N)^2 \frac {B^2}{Q^2N} \sum_{q\sim Q} \, \sum_{(a, q)=1} \Big|\sum_{k\in I} \mu(k)\frac 1{\phi(q)} \Big[\sum_{\Sb \Cal X\not\in\Cal B\\
(\mod q)\endSb} \overline{\Cal X}(a) c_{\Cal X}  (k)\Big]\Big|^2\tag 3.27
$$
where $I\subset[1,N]$ is an interval of size $\sim\frac{QN}B$.

By orthogonality
$$
(3.27) =(\log N)^2\frac {B^2}{Q^2N} \sum_{q\sim Q} \frac 1{\phi(q)} \sum_{\Sb \Cal X\not\in\Cal B\\ (\mod q)\endSb} \Big|\sum_{k\in I}
\mu(k) c_{\Cal X} (k)\Big|^2.\tag 3.28
$$
At this stage, we follow the argument from the proof of Theorem 3 in \cite{B-P}, noting that we have excluded characters induced by the
exceptional character $\Cal X^*$.

We briefly recall the steps.  Setting $d=(k, q)$, we have
$$
c_{\Cal X}(k) =\overline{\Cal X}_1 \Big(\frac kd\Big) c_{\Cal X} (d)\tag 3.29
$$
and $c_{\Cal X}(d)=0$ unless $q_1\Big|\frac qd, \frac q{dq_1}$ square free, $\Big(\frac q{q_1 d}, q_1\Big)=1$, in which case
$$
|c_{\Cal X}(d)| =\frac {\vp(q)}{\vp(\frac qd)} \, q_1^{1/2}.\tag 3.30
$$
Write
$$
\sum_{k\in I} \mu(k) c_{\Cal X}(k) =\sum_{d|\frac q{q_1}} c_{\Cal X} (d)\mu(d) \sum_{\Sb k_1\in \frac 1d I\\ (k_1, q)=1\endSb}
\mu(k_1)\overline{\Cal X}_1(k_1)
$$
hence
$$
\Big|\sum_{k\in I}\mu(k) c_{\Cal X} (k)\Big| \leq \sum_{d|\frac q{q_1}, d\, sf} \frac {\vp(q)}{\vp(\frac qd)} q_1^{1/2} \Big|\sum_{k_1\in \frac 1dI} \mu(k_1)\Cal X
(k_1)\Big|.\tag 3.31
$$
We use Theorem 4 from \cite{B-P} (the `case II-assumption' in the statement amounts to assuming $\Cal X_1\not= \Cal X^*$). Thus we have

\proclaim
{Lemma 3} Let $I_1\subset [1, N]$ be an interval of size $|I_1|>\frac NB$ and $\Cal X(\mod q), q\leq B$, not induced by the exceptional character.
Then
$$
\Big|\sum_{k\in I_1}\mu(k) \Cal X(k)\Big|\ll_A \frac {|I_1|}{(\log N)^A}.\tag 3.32
$$
\endproclaim

From (3.32)
$$
(3.31)\ll_A \frac {|I|}{(\log N)^A} \sum_{d|\frac q{q_1}, d \, sf} \frac {\vp(q)}{\vp(\frac qd)} \, q_1^{1/2} \, \frac 1d <(\log N)^{-A}\, 
q_1^{1/2} \, 2^{\omega(\frac q{q_1})} |I|.\tag 3.33
$$
Substituting (3.33) in (3.28) gives
$$
(3.28)\ll_A \frac 1{(\log N)^{A-2}} \, \frac BQ \, \sum_{q\sim Q} \, \frac 1{\vp(q)} \sum_{\Cal X(\mod q)} q_1^{\frac 12} \, 2^{\omega(\frac q{q_1})}
\Big|\sum_{k\in I} \mu(k)c_{\Cal X}(k)\Big|.\tag 3.34
$$
Using again (3.29), (3.30) and writing $q_2=\frac q{q_1}$,
$$
\align
(3.34)&< \frac 1{(\log N)^{A-2}} \, \frac BQ \, \sum_{q\sim Q} \, \frac 1{\vp(q)}\, \sum_{\Cal X(\mod q)} q_1^{\frac 12} \, 2^{\omega(\frac q{q_1})}
\sum_{d|\frac q{q_1}} |\mu(d)| \, |c_{\Cal X}(d)|\, \Big|\sum_{\Sb k_1\in \frac 1d I\\ (k_1, \frac q{q_1}) =1\endSb} \mu(k_1) \Cal X_1
(k_1)\Big|\\
&< \frac 1{(\log N)^{A-2}} \, \frac BQ \, \sum_{q_2\leq 2Q}\sum_{\Sb d|q_2\\ d \, sf\endSb}\frac {2^{\omega(q_2)}}{\vp(\frac {q_2} d)}
\sum_{q_1 \leq \frac {2Q}{q_2}} \frac {q_1}{\vp(q_1)} \underset{\Cal X(\mod q_1)}\to{\sum\nolimits^*} \Big|\sum_{\Sb k_1\in \frac 1dI\\ (k_1, q_2)=1\endSb}
\mu(k_1) \Cal X(k_1)\Big|.\tag 3.35
\endalign
$$
($\sum^*$ means summation over primitive characters).

According to \cite{B-P}, Theorem 5

\proclaim
{Lemma 4} Let $\ell\in\Bbb Z_+$. Then, with $R$ as in (3.2)
$$
\sum_{q\leq R}\, \frac q{\vp(q)}\,  \underset{\Cal X(\mod q)}\to{\sum\nolimits^*} \Big|\sum_{k\in I_1, (k, \ell)=1} \mu(k)\Cal X(k)\Big|\ll_\ve (\log N)^{74}
(R^2 N^{\frac 12}\ell^\ve+R|I_1|^{\frac 12} N^{\frac 3{10}} \ell^\ve +|I_1|).\tag 3.36
$$
\endproclaim

Hence
$$
\sum_{q_1\leq\frac {2Q}{q_2}} \frac {q_1}{\vp(q_1)} \underset {\Cal X(\mod q_1)}\to {\sum\nolimits^*} \Big|\sum_{k_1\in\frac 1d I, (k_1, q_2)=1} \mu (k_1)
\Cal X(k_1)\Big|\ll (\log N) ^{74} \frac {QN}{Bd}\tag 3.37
$$
and
$$
(3.26), (3.35) \ll_A (\log N)^{-A+76}N\sum_{q_2< 2Q} \, \frac {4^{\omega(q_2)}}{\vp(q_2)}\ll (\log N)^{-A+90} N.\tag 3.38
$$
This gives
$$
(3.24) \ll_A (\log N)^{-A}N.\tag 3.39
$$
Next, we analyze the contribution (3.25) of the exceptional character.

Estimate for some $Q<B$
$$
(3.25)\ll (\log N)^4 \, \frac B{NQ} \, \sum_I\sum_{q\sim Q, r|q} \frac 1{\vp(q)} \Big(\sum_{k\in I_1, k\, sf} |c_{\Cal X}(k)|\Big) \Big|\sum_{k\in I_2}
g(k) c_{\Cal X} (k)\Big|\tag 3.40
$$
where in the first sum, $I$ runs over a partition of $[1, N]$ in intervals of size $\sim \frac {QN}{B}$ and $I_1, I_2$ denote a pair of subintervals
of $I$.

Since $\Cal X$ is induced by $\Cal X^*(\mod r)$, it follows from (3.29), (3.30) that
$$
c_{\Cal X}(k) =\overline{\Cal X^*} \Big(\frac kd\Big) c_{\Cal X} (d) \ \text { and }\  |c_{\Cal X}(d)|<\frac {\vp(q)}{\vp(\frac qd)} r^{1/2}.\tag 3.41
$$
Hence
$$
\sum_{k\in I, k\, sf} |c_{\Cal X} (k) \leq r^{\frac 12} \sum_{k\in I, k\, sf} \Big(k, \frac qr\Big)
\ll r^{\frac 12}\, \frac {NQ} B \, 2^{\omega(\frac qr)}
\tag 3.42
$$
and
$$
(3.40)\ll (\log N)^4 r\sum_I\sum_{\Sb q\sim Q\\ r|q\endSb} 2^{\omega(\frac qr)} \sum_{d|\frac qr} \frac 1{\vp(\frac qd)} \Big|
\sum_{\Sb k\in I\\ d|k\endSb} g(k)\Cal X^*\Big(\frac kd\Big)\Big|.\tag 3.43
$$
Write $r=2^\gamma r_1, (r_1, 2)= 1, r_1>1$ ($r$ is a power of 2, cf. \cite{Gr}).
Let $\Cal X^*=\Cal X_0^*\Cal X_1^*, \Cal X_0^* (\mod 2^\gamma), \Cal X_1^* (\mod r_1)$.

Since $\Cal X_1^*$ is primitive,
$$
\Cal X_1^*(k)=\frac 1{\tau(\Cal X_1^*)} \sum_{(a, r_1)=1} \Cal X_1^*(a) e_{r_1} (ak)
$$
and we may also bound (3.43) by
$$
r\, r^{1/2}_1 (\log N)^4 \sum_I \sum_{\Sb q\sim Q\\ r|q\endSb} 2^{\omega (\frac qr)} \sum_{d|\frac qr} \frac 1{\vp(\frac qd)} \Big|\sum_{\Sb k\in I\\
d|k\endSb} g(k) \Cal X_0^* \Big(\frac kd\Big) e_{r_1} \Big(a\frac kd\Big)\Big|\tag 3.44
$$
with $(a, r_1)=1$.

Depending on whether $r_1$ is small or large, we use (3.44) or (3.43).

The estimates are carried out in the next section.

\bigskip
\noindent
{\bf 4. Major Arcs Analysis (II)}

First, observe that the `trivial' estimate on (3.43) gives
$$
(\log N)^5 \frac rQ\sum_{k\leq N} \sum_{d|k} d\, 2^{\omega(d)} \Big(\sum_{q_1\sim\frac Q{dr}} 2^{\omega(q_1)}\Big) |g(k)|.\tag 4.0
$$
Bounding
$$
\sum_{q_1\sim\frac Q{dr}} 2^{\omega(q_1)} \leq \Big(\frac Q{dr}\Big)^{1/2} \Big[\sum_{q_1 \sim \frac Q{dr}} 4^{\omega(q_1)}\Big]^{\frac 12}\ll
(\log N)^2 \frac Q{dr}
$$
we obtain
$$
\align
(4.0)& <(\log N)^7 \sum_{k<N} \Big(\sum_{d|k} 2^{\omega(d)}\Big) |g(k)|\\
&< (\log N)^7 \Big(\sum_{k<N} \tau (k)^4\Big)^{\frac 12}\Vert g\Vert_2\\
&<(\log N)^{24} \sqrt N \Vert g\Vert_2.\tag 4.1
\endalign
$$

Recall the definition (3.16) of $g$.

Let $J=J_\alpha =[n_1, n_2[\subset[0, n[$ where $n_2-n_1=m$ and $\frac n2> n_1>\frac n4$.

Write
$$
x=\sum_{j<n_1} x_j 2^j+\sum_{n_1\leq j< n_2} x_j 2^j +\sum_{n_2\leq j<n} x_j 2^j = u+2^{n_1} y+w.\tag 4.2
$$
and
$$
g(x)=\sum_{\Sb A\subset\{0, 1, \ldots, m-1\}\\ |A|<K_0\endSb} g_A(u, w) w_A(y)\tag 4.3
$$
where $y\in\{0, 1, \ldots, m-1\}$.

We will choose
$$
n_0 >\sqrt n \ \text { and } \ o(n)>m >\log B> C(\log n)n_0.\tag 4.4
$$
Recalling (3.13) and the choice $K\sim|\log c'|$, we have
$$
K_0\sim\Big(\frac 1{c'}\Big)^C \frac {n_0m}n\tag 4.5
$$
($C$ a constant).

Recalling (1.12), we need an approximation of $w_A$ with suitably restricted Fourier transform.
Let $\ell\in\Bbb Z_+$ be another parameter and replace the step function $h$ by a function $-1\leq h_0\leq 1$ satisfying
$$
h_0=\left\{\aligned 1+0(2^{-\ell}) \ & \text { if } \ 0\leq x<\frac 12-2^{-\ell}\\
-1+0(2^{-\ell}) \ & \text { if }  \ \frac 12 \leq x< 1\endaligned\right.
\tag 4.6
$$
and
$$
\supp \widehat h_0 \subset [-2^\ell, 2^\ell].\tag 4.7
$$
Denoting
$$
\tilde w_A (y) =\prod_{j\in A} h_0\Big(\frac y{2^{j+1}}\Big)\tag 4.8
$$
it follows from (4.6), (4.7) that
$$
2^{\frac {-m}2}\Vert w_A-\tilde w_A\Vert_2 < C|A|2^{-\ell/2}< CK_0 \, 2^{-\ell/2}\tag 4.9
$$
and
$$
\tilde w_A(y) =\sum_{\{ b_j\}_{j\in A}} \Big(\prod_{j\in A} \widehat h_0 (b_j)\Big) e\Big(\sum_{j\in A} \frac {b_j}{2^{j+1}} y\Big)\tag 4.10
$$
where
$$
\sum_{\{b_j\}}\Big|\prod_{j\in A}\widehat h_0 (b_j)\Big| \leq \Vert\widehat h_0\Vert_1^{|A|} \overset{(4.11)}\to < (c\ell)^{K_0}.\tag 4.11
$$
Note that if we replace each $w_A$ by $\tilde w_A$ in (4.3), we obtain a function $\tilde g$ that satisfies
$$
\align
\Vert g-\tilde g\Vert_{\ell^2[1, N]} &\leq \sum_A \Vert g_A\Vert_{\ell^2_{u, w}}\Vert w_A -\tilde w_A\Vert_{\ell^2_y}\\
&\overset{(4.9)}\to \ll K_0 2^{-\ell/2} \pmatrix m\\ K_0\endpmatrix \Vert g\Vert_2\\
&\lesssim n^{K_0} 2^{-\ell/2} \sqrt N.\tag 4.12
\endalign
$$

In view of the estimate (4.1), we can take
$$
\ell\sim(\log n)K_0.\tag 4.13
$$
Hence, by (4.7), (4.11)
$$
|b_j|< n^{CK_0} \ \text { and } \ \sum_{\{b_j\}} \Big|\prod_{j\in A} \widehat h_0(b_j)\Big| < \big(C(\log n)K_0\big)^{K_0}.\tag 4.14
$$
Returning to (3.43), (3.44), write
$$
d=2^\nu d_1; (d_1, 2)=1\tag 4.15
$$
and
$$
r=2^\gamma r_1; (r_1, 2)=1 \ \text { and } \ r_1>1.\tag 4.16
$$
We carry out the estimate by fixing $u, w$ and exploiting the $y$-variable.
By (4.2), since $\nu<\log Q=o(n)$, the condition $d|x$ clearly becomes
$$
2^\nu|u \ \text { and } \ y=z+d_1y', y'\in\Big\{0, 1, \ldots, \Big[\frac {2^m}{d_1}\Big]\Big\}\tag 4.17
$$
where
$$
z=z(u, w)\in \{0, 1, \ldots, d_1-1\} \ \text { is determined by \, $u+w+2^{n_1} z\equiv 0(\mod d_1)$}.
$$
Note also that $\log|I|>n-\log B>\big(1-o(1)\big)n$ and, by our choice of $n_1$, the restriction $x\in I$ essentially amounts to restricting
$w\in I$.

From the preceding and (4.3) with $w_A$ placed by $\tilde w_A$, we estimate in (3.43)
$$
\align
&\Big|\sum_{\Sb x\in I\\ d|x\endSb} g(x) \Cal X^* \Big(\frac xd\Big)\Big| \ \text { by }\\
&\sum_{w\in I} \, \sum_{2^\nu|u} \, \sum_{\Sb A\subset \{0, 1, \ldots, m-1\}\\ |A|<K_0\endSb} |g_A (u, w)|
\Big|\sum_{y'<\frac {2^m}{d_1}} \tilde w_A (z+d_1 y') \Cal X_1^* \Big(\frac {u+w+2^{n_1} z} d + 2^{n_1-\nu}y'\Big)\Big|\tag 4.18
\endalign
$$
and in (3.44)
$$
\align
&\Big|\sum_{\Sb x\in I\\ d|x\endSb} g(x) e_r \Big(\frac {ax}{d}\Big)\Big| \ \text { by }\\
&\sum_{w\in I} \, \sum_{2^\nu|u}\,  \sum_{A,|A| <K_0} |g_A(u, w)| \, \Big|\sum_{y'<\frac {2^m}{d_1}} \tilde w_A(z+d_1y') e\Big(\frac {a_1}{r_1}
y'\Big)\Big|\tag 4.19
\endalign
$$
where $a_1 \equiv 2^{n_1-\nu-\gamma} .a (\mod r_1), (a_1, r_1)=1$.

We consider first the case where $r_1$ is small, proceeding with (3.44), (4.19).
Using (4.10), the inner sum in (4.19) is bounded by
$$
\sum_{\{b_j\}}\prod_{j\in A} |\widehat h_0(b_j)| \ \Big|\sum_{y<\frac {2^m}{d_1}} e\Big(\Big(\frac {a_1}{r_1}+ d_1\sum_{J\in A} \frac
{b_j}{2^{j+1}}\Big)y\Big)\Big|.\tag 4.20
$$
Substitution in (3.44) gives then an estimate
$$
\align
&r_1^{\frac 12}\, r\, n^4 \frac 1Q\sum_A \sum_w\sum_\nu 2^\nu \sum_{2^\nu|u}\,  \sum_{d_1\lesssim 2^{-\nu} \frac Qr} 2^{\omega(d_1)} d_1
\Big(\sum_{q_1\sim \frac Q{dr}} 2^{\omega(q_1)}\Big) |g_A(u, \omega)|.(4.20)\\
&\ll r_1^{\frac 12} n^7 \sum_A\sum_w\sum_\nu\sum_{2^\nu|u}|g_A(u, w)|\sum_{d_1< 2^{-\nu}\frac Qr} 2^{\omega(d_1)}.(4.20)\\
&\overset{(4.10), (4.14)}\to \ll r_1^{\frac 12} n^{10} \pmatrix m\\ K_0\endpmatrix \big(C(\log n)K_0\big)^{K_0} \, 2^{-m}\Vert g\Vert_{\ell^1_x}.
(4.21) \\
&\ll r_1^{\frac 12} \big(Cn(\log n)K_0\big)^{K_0} N 2^{-m}.(4.21)\tag 4.22
\endalign
$$
(since $\Vert g\Vert_1 \leq \sqrt N\Vert g\Vert _2 < 2N)$

\noindent
and where (4.21) is an upper bound on
$$
\align
&\sum_{\Sb d\sim D\\ d \text { odd}\endSb } 2^{\omega (d)} \Big|\sum_{y<\frac {2^m}d} e\Big(\Big( \frac {a_1}{r_1} +d \sum_{j\in A}
\frac {b_j}{2^{j+1}}\Big)y\Big)\Big|\\
&\ll \sum_{\Sb d\sim D\\ d\text { odd}\endSb} 2^{\omega(d)} \Big[\Big\Vert \frac {a_1}{r_1} +d \Big( \sum_{j\in A} \frac {b_j}{2^{j+1}}
\Big)\Big\Vert +\frac D{2^m}\Big]^{-1}.\tag 4.23
\endalign
$$
We analyze (4.23) according to \cite {B2}, \S4.

Note that the contribution of those $d$ for which
$$
\Big\Vert\frac {a_1}{r_1}+d\Big(\sum_{j\in A} \frac{b_j}{2^{j+1}}\Big)\Big\Vert> B^2 n^{2K_0} 2^{-m}\tag 4.24
$$
in (4.22) is at most
$$
r_1^{\frac 12} \frac N{B^2} n^{-\frac 12K_0} \Big(\sum_{d\sim D} 2^{\omega (d)}\Big) < N.n^{-\frac 13 K_0}\tag 4.25
$$
and it remains to consider the others.

We assume \big(cf. (4.4)\big)
$$
m>10\log B+CK_0\log n.\tag 4.26
$$
Since $|A|<K_0$, there is some $m_1\in \big\{\big[\frac m4\big], \ldots, \big[\frac m2\big]\big\}$ such that
\hfill\break $A\cap [m_1, m_2]=\phi$ with $m_2=m_1+\big[\frac m{5K_0}]$.
Writing
$$
\beta =\sum_{j\in A} \frac{b_j}{2^{j+1}}=\sum_{\Sb j\in A\\ j<m_1\endSb} \frac {b_j}{2^{j+1}} +\sum_{\Sb j\in A\\j>m_2\endSb}
\frac {b_j}{2^{j+1}} =\beta_1+\beta_2
$$
it follows from (4.14) that
$$
|\beta_2|<n^{CK_0} 2^{-m_2}<n^{CK_0} 2^{-\frac m4}.\tag 4.27
$$
Hence, if $d$ fails (4.24), we have, since $(r_1, 2)=1$
$$
\frac 1{r_12^{m_1}} \leq \Big\Vert\frac {a_1}{r_1}+ d\beta_1\Vert < n^{CK_0} 2^{-m_2} d+2^{-\frac 34 m}.\tag 4.28
$$
by (4.25). Therefore
$$
d \geq  n^{-CK_0} \, 2^{[\frac m{5K_0}]} \frac 1{r_1} > 2^{\frac m{6K_0}}\tag 4.29
$$
provided
$$
r_1< 2^{\frac m{60K_0}}\tag 4.30
$$
and, recalling (4.5),
$$
m>CK_0^2 \log n \ \text { or } \ n> \big(\frac 1{c'}\big)^C(\log n) n_0\sqrt m.\tag 4.31
$$
We assume (4.30) and (4.31).

Hence, we may take $D> 2^{\frac m{6K_0}}$ in (4.23).

Next, denoting $d'$ the difference between distinct integers $d$ that fail (4.24), we have
$$
d'\not= 0, \Vert d'\beta\Vert< 2^{-\frac 34 m}.
$$
Hence
$$
\Vert d'\beta_1\Vert <|d'|. 2^{-m_1-\frac m{6K_0}}
$$
implying that either
$$
|d'|> 2^{\frac m{6K_0}}\tag 4.32
$$
or
$$
0<|d'|\leq 2^{\frac m{6K_0}} \ \text { and } \  d'\beta_1 \equiv 0 (\mod 1).\tag 4.33
$$

Let $\tau$ be the power of 2 in the factorization of $d'$.
Then 
$$
2^\tau \leq 2^{\frac m{6K_0}} \ \text {  and } \ 2^\tau \beta_1 \equiv 0 (\mod 1).\tag 4.34
$$
If $d\sim D$ is any integer failing (4.24), it follows from (4.34) that
$$
\frac 1B<\frac 1{r_1} \leq\Big\Vert\frac {2^\tau a_1}{r_1}\Big\Vert\leq 2^\tau\Big\Vert\frac {a_1}{r_1} +
 d\beta\Big\Vert+ 2^\tau D|\beta_2| < 2^{-\frac m2}
+2^{\frac m{6K_0}} Bn^{CK_0} 2^{- \frac m4} <  2^{-\frac m8}\tag 4.35
$$
(contradiction).

This shows that the difference between two distinct elements $d$ failing (4.24) is at least $2^{\frac m{6K_0}}$.
Therefore, since $D> 2^{\frac m{6K_0}}$, we get
$$
\align
&\sum_{\Sb d\sim D\\ d  \text { fails (4.24)}\endSb} 2^{\omega(d)} \frac 1{\Vert \frac {a_1}{r_1} +d\beta\Vert +\frac D{2^m}}\leq\\
&\frac {2^m} D\Big(\sum_{d\sim D} 4^{\omega(d)}\Big)^{\frac 12} |\{ d\sim D; d \text { fails (4.24)}|^{\frac 12}\lesssim\\
&\frac {2^m}D(\log D)^4 D^{\frac 12} \Big(1+D. 2^{-\frac m{6K_0}}\Big)^{\frac 12}\lesssim\\
&2^m(\log D)^4 \, 2^{-\frac m{12K_0}}.\tag 4.36
\endalign
$$
Substituting (4.36) in (4.22) gives the estimate
$$
r_1^{\frac 12} (Cn(\log n)K_0)^{K_0+4} 2^{-\frac m{12K}} N\overset {(4.26)}\to < 2^{-\frac m{13K_0}} r_1^{\frac 12} N.\tag 4.37
$$
Summarizing, from (4.25), (4.37), (4.5), and by assumption (4.30) on $r$, we have 
$$
(3.44)< \big(n^{-\frac 13K_0}+r_1^{\frac 12} 2^{-\frac m{13K_0}}\big) N< 2n^{-\frac 13K_0}N< 2^{-\frac{n_0m}n\log n} N.\tag 4.38
$$
Recall also conditions (2.6), (4.26), (4.31) on $B$ and $m$
$$
\log B> C(\log n)n_0 \,  m> 10\log B \ \text { and } \ n>\Big(\frac 1{c'}\Big)^C n_0\sqrt m\log n.\tag 4.39
$$
Next, consider the case that $r_1$ is `large', in the sense that (4.30) fails
$$
r_1\geq 2^{\frac m{60K_0}}.\tag 4.40
$$
We now proceed with (3.43), (4.18), obtaining instead of (4.22) the bound
$$
\big(Cn(\log n)K_0\big)^{K_0} N2^{-m}.(4.41)\tag 4.42
$$
with (4.41) an upper bound on
$$
\sum_{\Sb d\sim D\\ d\text { odd}\endSb} 2^{\omega(d)}\Big|\sum_{y<\frac {2^m}d} e\Big(d\Big(\sum_{j\in A} \frac {b_j}{2^{j+1}}\Big)y\Big) \Cal X_1^*
	(\ell+2^ty)\Big|.\tag 4.41
$$
for some $\ell, t\in\Bbb Z_+$.

The inner sum in (4.41) may be further evaluated by
$$
r_1+\frac {2^m}{dr_1}.(4.43)\tag 4.42
$$
with (4.43) of the form
$$
\Big|\sum_{y=0}^{r_1-1} e(\theta y)\Cal X_1^* (\ell_1+2^ty)\Big|=(4.43)
$$
for some $\theta \in\Bbb R$ and $\ell_1\in\Bbb Z$.

Writing $\theta =\frac b{r_1} +\psi(\mod 1), |\psi|<\frac 1{r_1}$, partial summation gives an estimate on (4.43) by an incomplete sum
$$
\sum_{y\in V} e_{r_1} (by)\Cal X_1^* (\ell_1+2^t y)\tag 4.44
$$
with $V$ a subinterval of $\{0, 1, \ldots, r_1-1\}$.

Completing the sum (4.44), we obtain, since $(r_1, 2)=1$
$$
\align
&(\log r_1)\Big|\sum_{y=0}^{r_1-1} e_{r_1} (b_1 y) \Cal X_1^* (\ell_1+ 2^ty)\Big|=\\
&(\log r_1)\Big|\sum^{r_1-1}_{y=0} e_{r_1} (b_2y)\Cal X_1^* (y)\Big| =(\log r_1)\sqrt {r_1}
\endalign
$$
since $\Cal X_1^*$ is primitive.

Therefore (4.41) is bounded by
$$
\sum_{d\sim D} 2^{\omega(d)} \Big(r_1+\frac {2^m}{d\sqrt{r_1}}\log r_1\Big) \ll (\log N)^4 \Big(B^2+ \frac {2^m}{\sqrt{r_1}}\Big)
\overset {(4.39), (4.40)}\to< 2^m 2^{-\frac m{121K_0}}
$$
and, by (4.31)
$$
(4.42)< N 2^{-\frac m{122K_0}}.\tag 4.45
$$
Thus the estimate (4.38) also holds for $r_1$ large.

In view of (4.39), choose $m, \log B\sim n_0\log n$ and $n_0$ such that
$$
n^{1/2} < n_0<(c')^C\frac {n^{2/3}}{\log n}.\tag 4.46
$$
Since (4.38) is an estimate for (3.25) and recalling (3.39), this proves that
$$
\Big|\int_{\Cal M} S(\alpha) \overline S_g(\alpha) d\alpha\Big| \ll_A \frac N{(\log N)^A}.\tag 4.47
$$
Recalling (3.19) and the construction from \S1, we therefore proved that
$$
\sum_{|A|\leq n_0}|\widehat\mu(A)|^2 <(\log n)^{-\lambda} \ \text { provided} \ n_0<(\log n)^{-C\lambda} \, n^{2/3}\tag 4.48
$$
for any fixed constant $\lambda\geq 1$.

\enddocument